\documentclass[12pt]{article}
\usepackage{amsfonts}

\newtheorem{proposition}{Proposition}

\begin{document}

\renewcommand{\theequation}{\thesection.\arabic{equation}}

\setcounter{equation}{0}
\title{The Weierstrass representation of discrete isotropic surfaces in $\mathbb{R}^{2,1}$, $\mathbb{R}^{3,1}$
and $\mathbb{R}^{2,2}$}

\author{D.V.Zakharov \thanks{Columbia University, New York, USA; e-mail:
zakharov@math.columbia.edu}}

\date{}

\maketitle
\begin{abstract}
Using an integrable discrete Dirac operator, we construct a discrete version of the Weierstrass representation of time-like
surfaces parametrized along isotropic directions in $\mathbb{R}^{2,1}$, $\mathbb{R}^{3,1}$ and $\mathbb{R}^{2,2}$.
The corresponding discrete surfaces have isotropic edges. We show that any discrete surface satisfying a general monotonicity
condition and having isotropic edges admits such a representation.

\end{abstract}

\section{Introduction}

\setcounter{equation}{0}

The classical Weierstrass representation \cite{taimanov2}, \cite{konopelchenko1},
\cite{konopelchenko-taimanov} associates to every solution of the Dirac equation
\begin{equation}
\partial_z \varphi=u\psi,\quad \partial_{\bar{z}}\psi=-u\varphi,\quad u=\bar{u}
\label{Dirac}
\end{equation}
a conformally embedded surface in $\mathbb{R}^3$, and any conformally embedded surface locally admits such
a description. The mean curvature of a surface embedded under the Weierstrass representation admits a particularly simple description in terms of the Dirac operator (\ref{Dirac}), which makes the Weierstrass representation a powerful tool for
studying surfaces of constant mean curvature and minimal surfaces. In addition, the Weierstrass representation has received significant attention in recent years as a possible approach to the Willmore conjecture (see the recent survey
\cite{taimanov3}). 

Many different versions of the Weierstrass representation have been found. In particular, various reductions
of the generalized Dirac operator 
\begin{equation}
\partial_{\xi}\varphi=u\psi,\quad \partial_{\eta}\psi=v\varphi,
\label{gen-Dirac}
\end{equation}
where $\xi$ and $\eta$ are complex variables and $u$ and $v$ are complex-valued functions, can be used to
construct surfaces in three- and four-dimensional Euclidean and pseudo-Euclidean spaces
(see \cite{konopelchenko2}, \cite{konopelchenko3}, \cite{konopelchenko4}). There are essentially
two different reductions. The reduction $\eta=\bar{\xi}$ describes conformally embedded {\it space-like} surfaces in  $\mathbb{R}^{4}$, $\mathbb{R}^{3,1}$ and $\mathbb{R}^{2,2}$, while the reduction $\eta=\bar{\eta}$, $\xi=\bar{\xi}$ describes  describes {\it time-like} surfaces in
$\mathbb{R}^{2,1}$, $\mathbb{R}^{3,1}$ and $\mathbb{R}^{2,2}$ with isotropic coordinate lines. 

In \cite{zakharov}, the author considered the following discrete operator as an integrable discretization of the generalized Dirac operator (\ref{gen-Dirac}):
\begin{equation}
\tau_2\varphi=\alpha\varphi+\beta\psi,\quad \tau_1\psi=\gamma\varphi+\delta\psi.
\end{equation}
Here the functions $\varphi$ and $\psi$ depend on two discrete variables and $\tau_1$ and $\tau_2$ denote the
translation operators in these variables. It is natural to ask whether this discretization can
be used to construct a discrete analogue of the Weierstrass representation.

In this paper, we construct a discretization of the Weierstrass representation of time-like surfaces
with isotropic coordinate lines. We show that solutions of a certain discrete Dirac operator can be used to construct lattices in $\mathbb{R}^{2,1}$, $\mathbb{R}^{3,1}$ and $\mathbb{R}^{2,2}$ with the geometric property that
every edge of the lattice is an isotropic vector. We also show that any such lattice satisfying a
generic condition of monotonicity can be described in this way. 

As mentioned above, the Weierstrass representation of conformal surfaces in Euclidean and pseudo-Euclidean
spaces involves a different Dirac operator, involving $\partial$ and $\bar{\partial}$ instead of two real variables.
The direct approach taken in this paper does not work in this case. For a general reference on
discrete differential geometry and integrable lattices, see \cite{bobenko-suris}.

The plan of the paper is as follows. 
In Section 2, we recall the Weierstrass representation of time-like surfaces in $\mathbb{R}^{2,1}$,
$\mathbb{R}^{3,1}$ and $\mathbb{R}^{2,2}$, following \cite{konopelchenko2}. In Section 3, we show that for each of the three cases, an appropriate discrete analogue of the Dirac operator can be used
to construct, using the same formulas as in the differential case,  discrete surfaces in $\mathbb{R}^{2,1}$, $\mathbb{R}^{3,1}$ and $\mathbb{R}^{2,2}$ with isotropic edges, and we show that any such surface can be described in this way. In Section 4, we show that by letting the mesh size tend to zero we obtain the continuous construction.

\section{The Weierstrass representation of time-like isotropic surfaces}

\setcounter{equation}{0}

In this section, we recall the Weierstrass representation of time-like surfaces in $\mathbb{R}^{2,1}$, $\mathbb{R}^{3,1}$ and
$\mathbb{R}^{2,2}$.

Let $\mathbb{R}^{n,m}$ denote pseudo-Euclidean space of dimension $n+m$ with a metric $\langle\cdot,\cdot\rangle$ of signature $(n,m)$, and let 
$S$ be a surface. We say that an embedding of $S$ in $\mathbb{R}^{n,m}$ is {\it time-like} if the induced metric on $S$ has signature $(1,1)$. In the case when $(n,m)=(2,1)$, $(3,1)$ or $(2,2)$, any such embedding can be locally described
using solutions of various reductions of the following Dirac equation:
\begin{equation}
\partial_y \varphi=p \psi,\quad
\partial_x \psi=q\varphi,
\label{realdirac}
\end{equation}
where $x$ and $y$ are real variables. We describe these three cases individually.

\subsection{The $\mathbb{R}^{2,1}$-case}

Suppose that the complex-valued functions $\varphi$, $\psi$ are defined on some simply-connected domain
$U\subset\mathbb{R}^2$ and satisfy the equations
\begin{equation}
\partial_y \varphi=p \psi,\quad\partial_x \psi=p\varphi, \quad \bar{p}=p.
\label{cdirac21}
\end{equation}
Let $P$ be a point in $U$. Then the formulas
\begin{equation}
X^1(Q)=\frac{1}{2}\displaystyle\int_P^Q\left[(\varphi^2+\bar{\varphi}^2)d x+(\psi^2+\bar{\psi}^2)d y\right],
\end{equation}
\begin{equation}
X^2(Q)=\frac{i}{2}\displaystyle\int_P^Q\left[(\varphi^2-\bar{\varphi}^2)d x+(\psi^2-\bar{\psi}^2)d y\right],
\end{equation}
\begin{equation}
X^3(Q)=\displaystyle\int_P^Q\left[\varphi\bar{\varphi}d x+\psi\bar{\psi}d y\right],
\end{equation}
define an embedding $\vec{X}:U\rightarrow\mathbb{R}^{2,1}$, such that the induced metric on $U$ has signature $(+,-)$ and moreover the directions $x=$ const, $y=$ const are isotropic, i.e.
\begin{equation}
\left\langle \frac{\partial \vec{X}}{\partial x},\frac{\partial \vec{X}}{\partial x}\right\rangle=0,\quad
\left\langle \frac{\partial \vec{X}}{\partial y},\frac{\partial \vec{X}}{\partial y}\right\rangle=0.
\end{equation}
Conversely, any time-like surface embedded in $\mathbb{R}^{2,1}$ locally
admits such a representation.

\subsection{The $\mathbb{R}^{3,1}$-case}

Suppose that the complex-valued functions $\varphi_i$, $\psi_i$, $i=1,2$ are defined on some simply-connected domain
$U\subset\mathbb{R}^2$ and satisfy the equations
\begin{equation}
\partial_y \varphi_i=p \psi_i,\quad
\partial_x \psi_i=\bar{p}\varphi_i.
\label{cdirac31}
\end{equation}
Let $P$ be a point in $U$. Then the formulas
\begin{equation}
X^1(Q)=\frac{1}{2}\displaystyle\int_P^Q\left[(\varphi_1\bar{\varphi}_2+\bar{\varphi}_1\varphi_2)d x+
(\psi_1\bar{\psi}_2+\bar{\psi}_1\psi_2)d y\right],
\end{equation}
\begin{equation}
X^2(Q)=\frac{i}{2}\displaystyle\int_P^Q\left[(\varphi_1\bar{\varphi}_2-\bar{\varphi}_1\varphi_2)d x+
(\psi_1\bar{\psi}_2-\bar{\psi}_1\psi_2)d y\right],
\end{equation}
\begin{equation}
X^3(Q)=\frac{1}{2}\displaystyle\int_P^Q\left[(\varphi_1\bar{\varphi}_1-\varphi_2\bar{\varphi}_2)d x+
(\psi_1\bar{\psi}_1-\psi_2\bar{\psi}_2)d y\right],
\end{equation}
\begin{equation}
X^4(Q)=\frac{1}{2}\displaystyle\int_P^Q\left[(\varphi_1\bar{\varphi}_1+\varphi_2\bar{\varphi}_2)d x+
(\psi_1\bar{\psi}_1+\psi_2\bar{\psi}_2)d y\right],
\end{equation}
define an embedding $\vec{X}:U\rightarrow\mathbb{R}^{3,1}$, such that the induced metric on $U$ has signature $(+,-)$ and moreover the directions $x=$ const, $y=$ const are isotropic. Conversely, any time-like surface embedded in $\mathbb{R}^{3,1}$ locally
admits such a representation.

\subsection{The $\mathbb{R}^{2,2}$-case}

Suppose that the functions $\varphi_i$, $\psi_i$, $i=1,2$
are defined on some simply-connected domain $U\subset\mathbb{R}^2$ and satisfy the equations
\begin{equation}
\partial_y \varphi_1=p\psi_1,\quad \partial_x \psi_1=q\varphi_1,\quad
\partial_y\varphi_2=q\psi_2,\quad \partial_x\psi_2=p\varphi_2,\quad \bar{p}=p,\quad \bar{q}=q.
\label{cdirac22}
\end{equation}
Let $P$ be a point in $U$. Then the formulas
\begin{equation}
X^1(Q)=\frac{1}{2}\displaystyle\int_P^Q\left[(\varphi_1\varphi_2+\bar{\varphi}_1\bar{\varphi}_2)d x+
(\psi_1\psi_2+\bar{\psi}_1\bar{\psi}_2)d y\right],
\end{equation}
\begin{equation}
X^2(Q)=\frac{i}{2}\displaystyle\int_P^Q\left[(\varphi_1\varphi_2-\bar{\varphi}_1\bar{\varphi}_2)d x+
(\psi_1\psi_2-\bar{\psi}_1\bar{\psi}_2)d y\right],
\end{equation}
\begin{equation}
X^3(Q)=\frac{1}{2}\displaystyle\int_P^Q\left[(\varphi_1\bar{\varphi}_2+\bar{\varphi}_1\varphi_2)d x+
(\psi_1\bar{\psi}_2+\bar{\psi}_1\psi_2d y\right],
\end{equation}
\begin{equation}
X^4(Q)=\frac{i}{2}\displaystyle\int_P^Q\left[(\varphi_1\bar{\varphi}_2-\bar{\varphi}_1\varphi_2)d x+
(\psi_1\bar{\psi}_2-\bar{\psi}_1\psi_2d y\right],
\end{equation}
define an embedding $\vec{X}:U\rightarrow\mathbb{R}^{2,2}$, such that the induced metric on $U$ has signature $(+,-)$ and moreover the directions $x=$ const, $y=$ const are isotropic. Conversely, any time-like surface embedded in $\mathbb{R}^{2,2}$ locally
admits such a representation.

\section{Discrete isotropic surfaces}

\setcounter{equation}{0}

As we saw in the previous section, solutions of certain reductions of the Dirac equation (\ref{realdirac})
can be used to construct time-like embeddings of $\mathbb{R}^2$ into $\mathbb{R}^{2,1}$, $\mathbb{R}^{3,1}$ and
$\mathbb{R}^{2,2}$, with the property that all coordinate lines are isotropic. In the paper \cite{zakharov},
the author considered a discretization of the spectral data associated to the generalized Dirac operator (\ref{gen-Dirac})
and derived the following discrete system, which can be considered as a generalized discrete Dirac operator:
\begin{equation}
\tau_2\varphi=\alpha\varphi+\beta\psi,\quad \tau_1\psi=\gamma\varphi+\delta\psi.
\label{gen-discrete-Dirac}
\end{equation}
Here $\varphi$ and $\psi$ are functions of two discrete variables $(n,m)\in\mathbb{Z}^2$ and $\tau_1$ and $\tau_2$ are the
translation operators in these variables. By considering symmetries on the spectral data similar to the ones which
reduce the generalized operator (\ref{gen-Dirac}) to the standard one (\ref{Dirac}), the author obtained the following
reduction of (\ref{gen-discrete-Dirac}):
\begin{equation}
\tau_2\varphi=\alpha\varphi+\beta\psi,\quad \tau_1\psi=\beta\varphi+\alpha\psi,\quad\alpha^2-\beta^2=1.
\label{discrete-Dirac}
\end{equation}
In this section we show how that certain reductions of the generalized discrete Dirac operator (\ref{gen-discrete-Dirac})
similar to the one above can be used to construct, using the same formulas as in the differential
case, discrete surfaces in $\mathbb{R}^{2,1}$, $\mathbb{R}^{3,1}$ and $\mathbb{R}^{2,2}$ with isotropic edges.

A {\it discrete surface} in $\mathbb{R}^{n,m}$ is a map $\vec{X}:\mathbb{Z}^2\rightarrow\mathbb{R}^{n,m}$. The {\it edges} of 
a discrete surface $\vec{X}$ are the vectors
\begin{equation}
\vec{F}(n,m)=\vec{X}(n+1,m)-\vec{X}(n,m),\quad\vec{G}(n,m)=\vec{X}(n,m+1)-\vec{X}(n,m).
\end{equation}
Conversely, a pair of functions $\vec{F}:\mathbb{Z}^2\rightarrow\mathbb{R}^{n,m}$ and
$\vec{G}:\mathbb{Z}^2\rightarrow\mathbb{R}^{n,m}$ defines a discrete surface (up to translation) if and only if it satisfies the
following consistency condition:
\begin{equation}
\vec{F}(n,m+1)-\vec{F}(n,m)=\vec{G}(n+1,m)-\vec{G}(n,m).
\label{consistency}
\end{equation}
We use the following notation for surfaces defined in terms of their edges:
\begin{equation}
\vec{X}=\displaystyle\sum\left(\vec{F}\Delta_1+\vec{G}\Delta_2\right).
\end{equation}
A discrete surface is called {\it non-degenerate} if its edges are linearly independent at every lattice point. A discrete
surface in $\mathbb{R}^{n,m}$ is called {\it isotropic} if all of its edges are light-like vectors:
\begin{equation}
\langle \vec{F}(n,m),\vec{F}(n,m)\rangle=0,\quad\langle \vec{G}(n,m),\vec{G}(n,m)\rangle=0.
\label{isotropic}
\end{equation}

The main result of this paper is that isotropic discrete surfaces in $\mathbb{R}^{2,1}$, $\mathbb{R}^{3,1}$ and
$\mathbb{R}^{2,2}$ satisfying a certain monotonicity condition are described by solutions of a discrete Dirac equation.

\subsection{The $\mathbb{R}^{2,1}$-case}

\begin{proposition}
Suppose that the functions $\varphi$, $\psi$ satisfy the following discrete Dirac equation:
\begin{equation}
\tau_2\varphi=\alpha\varphi+\beta\psi,\quad \tau_1\psi=\beta\varphi+\alpha\psi, \quad \bar{\alpha}=\alpha,\quad\bar{\beta}=\beta,
\quad \alpha^2-\beta^2=1.
\label{dirac21}
\end{equation}
Then the formulas
\begin{equation}
X_1=\frac{1}{2}\displaystyle\sum\left[(\varphi^2+\bar{\varphi}^2)\Delta_1+(\psi^2+\bar{\psi}^2)\Delta_2\right],
\label{sum21a}
\end{equation}
\begin{equation}
X_2=\frac{i}{2}\displaystyle\sum\left[(\varphi^2-\bar{\varphi}^2)\Delta_1+(\psi^2-\bar{\psi}^2)\Delta_2\right],
\label{sum21b}
\end{equation}
\begin{equation}
X_3=\displaystyle\sum\left[\varphi\bar{\varphi}\Delta_1+\psi\bar{\psi}\Delta_2\right],
\label{sum21c}
\end{equation}
define an isotropic discrete surface $\vec{X}:\mathbb{Z}^2\rightarrow\mathbb{R}^{2,1}$. Conversely, if
$\vec{X}:\mathbb{Z}^2\rightarrow\mathbb{R}^{2,1}$ is a non-degenerate isotropic discrete surface that satisfies the following
condition
\begin{equation}
X_3(n+1,m)-X_3(n,m)>0,\quad X_3(n,m+1)-X_3(n,m)>0\mbox{ for all }(n,m)\in\mathbb{Z}^2,
\label{monotonic21}
\end{equation}
then there exist functions $\varphi$ and $\psi$ satisfying equation (\ref{dirac21}) such that equations (\ref{sum21a})-(\ref{sum21c}) hold.
\end{proposition}

\noindent {\bf Proof.} Given functions $\varphi$, $\psi$ satisfying (\ref{dirac21}), a direct calculation
shows that the edges given by equations (\ref{sum21a})-(\ref{sum21c}) are isotropic (\ref{isotropic}) and satisfy
the consistency condition (\ref{consistency}), hence define an isotropic discrete surface in $\mathbb{R}^{2,1}$.

Conversely, suppose that $\vec{X}:\mathbb{Z}^2\rightarrow\mathbb{R}^{2,1}$ is an isotropic discrete surface satisfying the
monotonic condition (\ref{monotonic21}). The edges of the lattice satisfy the equations
\begin{equation}
F_1^2+F_2^2=F_3^2,\quad G_1^2+G_2^2=G_3^2,\quad F_3>0,\quad G_3>0,
\end{equation}
therefore, there exist functions $\varphi$ and $\psi$, defined up to multiplication by $\pm 1$, such that the edges are given by the formulas (\ref{sum21a})-(\ref{sum21c}). The consistency condition (\ref{consistency}) implies that these functions satisfy the following equations
\begin{equation}
(\tau_2\varphi)^2-\varphi^2=(\tau_1\psi)^2-\psi^2,
\label{consistency21a}
\end{equation}
\begin{equation}
(\tau_2\varphi)(\tau_2\bar{\varphi})-\varphi\bar{\varphi}=(\tau_1\psi)(\tau_1\bar{\psi})-\psi\bar{\psi},
\label{consistency21b}
\end{equation}
and the non-degeneracy condition implies that
\begin{equation}
\varphi\bar{\psi}-\bar{\varphi}\psi\neq 0.
\label{nondeg21}
\end{equation}

The above equation implies that there exist unique real-valued functions $\alpha$, $\beta$, $\gamma$ and $\delta$ such
that the following system of equations is satisfied:
\begin{equation}
\tau_2\varphi=\alpha\varphi+\beta\psi,\quad
\tau_1\psi=\gamma\varphi+\delta\psi.
\end{equation}
Solving this system, we get
\begin{equation}
\alpha=\frac{\bar{\psi}(\tau_2\varphi)-\psi(\tau_2\bar{\varphi})}{\varphi\bar{\psi}-\bar{\varphi}\psi},\quad
\beta=\frac{\varphi(\tau_2\bar{\varphi})-\bar{\varphi}(\tau_2\varphi)}{\varphi\bar{\psi}-\bar{\varphi}\psi},
\end{equation}
\begin{equation}
\gamma=\frac{\bar{\psi}(\tau_1\psi)-\psi(\tau_1\bar{\psi})}{\varphi\bar{\psi}-\bar{\varphi}\psi},\quad
\delta=\frac{\varphi(\tau_1\bar{\psi})-\bar{\varphi}(\tau_1\psi)}{\varphi\bar{\psi}-\bar{\varphi}\psi},
\end{equation}
and a direct calculation using (\ref{consistency21a})-(\ref{consistency21b}) shows that
\begin{equation}
\alpha^2-\gamma^2=1,\quad\delta^2-\beta^2=1,\quad\alpha\beta=\gamma\delta.
\end{equation}
Solving this system we get that $\delta=\lambda\alpha$ and $\gamma=\lambda\beta$, where $\lambda=\pm 1$. Changing the
signs of $\psi$ at every point if necessary, we can set $\lambda=1$, so that the functions $\varphi$ and $\psi$ satisfy
the system (\ref{dirac21}). This proves the proposition.

\subsection{The $\mathbb{R}^{3,1}$-case}

\begin{proposition}

Suppose that the functions $\varphi_i$, $\psi_i$, $i=1,2$ satisfy the following discrete Dirac
equation:
\begin{equation}
\tau_2\varphi_i=\alpha\varphi_i+\beta\psi_i,\quad
\tau_1\psi_i=\bar{\beta}\varphi_i+\bar{\alpha}\psi_i,\quad |\alpha|^2-|\beta|^2=1.
\label{dirac31}
\end{equation}
Then the formulas
\begin{equation}
X_1=\frac{1}{2}\displaystyle\sum\left[(\varphi_1\bar{\varphi}_2+\bar{\varphi}_1\varphi_2)\Delta_1+
(\psi_1\bar{\psi}_2+\bar{\psi}_1\psi_2)\Delta_2\right],
\label{sum31a}
\end{equation}
\begin{equation}
X_2=\frac{i}{2}\displaystyle\sum\left[(\varphi_1\bar{\varphi}_2-\bar{\varphi}_1\varphi_2)\Delta_1+
(\psi_1\bar{\psi}_2-\bar{\psi}_1\psi_2)\Delta_2\right],
\label{sum31b}
\end{equation}
\begin{equation}
X_3=\frac{1}{2}\displaystyle\sum\left[(\varphi_1\bar{\varphi}_1-\varphi_2\bar{\varphi}_2)\Delta_1+
(\psi_1\bar{\psi}_1-\psi_2\bar{\psi}_2)\Delta_2\right],
\label{sum31c}
\end{equation}
\begin{equation}
X_4=\frac{1}{2}\displaystyle\sum\left[(\varphi_1\bar{\varphi}_1+\varphi_2\bar{\varphi}_2)\Delta_1+
(\psi_1\bar{\psi}_1+\psi_2\bar{\psi}_2)\Delta_2\right],
\label{sum31d}
\end{equation}
define an isotropic discrete surface $\vec{X}:\mathbb{Z}^2\rightarrow\mathbb{R}^{3,1}$. Conversely, if
$\vec{X}:\mathbb{Z}^2\rightarrow\mathbb{R}^{3,1}$ is a non-degenerate isotropic discrete surface that satisfies the following
condition
\begin{equation}
X_4(n+1,m)-X_4(n,m)>0,\quad X_4(n,m+1)-X_4(n,m)>0\mbox{ for all }(n,m)\in\mathbb{Z}^2,
\label{monotonic31}
\end{equation}
then there exist functions $\varphi_i$, $\psi_i$, $i=1,2$ satisfying equation (\ref{dirac31})such that equations (\ref{sum31a})-(\ref{sum31d}) hold.
\end{proposition}

\noindent {\bf Proof.} Given functions $\varphi_i$, $\psi_i$, $i=1,2$ satisfying (\ref{dirac31}), a direct calculation
shows that the edges given by equations (\ref{sum31a})-(\ref{sum31d}) are isotropic (\ref{isotropic}) and satisfy
the consistency condition (\ref{consistency}), hence define an isotropic discrete surface in $\mathbb{R}^{3,1}$.

Conversely, suppose that $\vec{X}:\mathbb{Z}^2\rightarrow\mathbb{R}^{3,1}$ is an isotropic discrete surface satisfying the
monotonic condition (\ref{monotonic31}). The edges of the lattice satisfy the equations
\begin{equation}
F_1^2+F_2^2+F_3^2=F_4^2,\quad G_1^2+G_2^2+G_3^2=G_4^2,\quad F_4>0,\quad G_4>0,
\end{equation}
therefore, there exist functions $\varphi_i$ and $\psi_i$, where $i=1,2$, such that the edges
are given by the formulas (\ref{sum31a})-(\ref{sum31d}). These functions are defined up to the following local gauge equivalence:
\begin{equation}
\varphi_1\rightarrow e^{i\zeta}\varphi_1,\quad\varphi_2\rightarrow e^{i\zeta}\varphi_2,\quad
\psi_1\rightarrow e^{i\xi}\psi_1,\quad\psi_2\rightarrow e^{i\xi}\varphi_2,
\label{gauge31}
\end{equation}
where $\zeta$ and $\xi$ are real-valued functions. The consistency condition (\ref{consistency}) implies that these functions satisfy the following equations
\begin{equation}
(\tau_2\varphi_1)(\tau_2\bar{\varphi}_1)-\varphi_1\bar{\varphi}_1=(\tau_1\psi_1)(\tau_1\bar{\psi}_1)-\psi_1\bar{\psi}_1,
\label{consistency31a}
\end{equation}
\begin{equation}
(\tau_2\varphi_1)(\tau_2\bar{\varphi}_2)-\varphi_1\bar{\varphi}_2=(\tau_1\psi_1)(\tau_1\bar{\psi}_2)-\psi_1\bar{\psi}_2,
\label{consistency31b}
\end{equation}
\begin{equation}
(\tau_2\varphi_2)(\tau_2\bar{\varphi}_2)-\varphi_2\bar{\varphi}_2=(\tau_1\psi_2)(\tau_1\bar{\psi}_2)-\psi_2\bar{\psi}_2,
\label{consistency31c}
\end{equation}
and the non-degeneracy condition implies that
\begin{equation}
\varphi_1\psi_2-\varphi_2\psi_1\neq 0.
\end{equation}

The above equation implies that there exist unique complex-valued functions $\alpha$, $\beta$, $\gamma$ and $\delta$ such
that the following system of equations is satisfied:
\begin{equation}
\tau_2\varphi_i=\alpha\varphi_i+\beta\psi_i,\quad \tau_1\psi_i=\gamma\varphi_i+\delta\psi_i,\quad i=1,2.
\end{equation}
We can explicitly solve these to obtain
\begin{equation}
\alpha=\frac{\psi_2(\tau_2\varphi_1)-\psi_1(\tau_2\varphi_2)}{\varphi_1\psi_2-\varphi_2\psi_1},\quad
\beta=\frac{\varphi_1(\tau_2\varphi_2)-\varphi_2(\tau_2\varphi_1)}{\varphi_1\psi_2-\varphi_2\psi_1},
\end{equation}
\begin{equation}
\gamma=\frac{\psi_2(\tau_1\psi_1)-\psi_1(\tau_1\psi_2)}{\varphi_1\psi_2-\varphi_2\psi_1},\quad
\delta=\frac{\varphi_1(\tau_1\psi_2)-\varphi_2(\tau_1\psi_1)}{\varphi_1\psi_2-\varphi_2\psi_1}.
\end{equation}
and a direct calculation using (\ref{consistency31a})-(\ref{consistency31c}) shows that
\begin{equation}
\alpha\bar{\alpha}-\gamma\bar{\gamma}=1,\quad\delta\bar{\delta}-\beta\bar{\beta}=1,\quad
\alpha\bar{\beta}-\gamma\bar{\delta}=0.
\end{equation}
Solving this system we get that $\delta=\lambda\bar{\alpha}$ and $\gamma=\lambda\bar{\beta}$, where
$\lambda\bar{\lambda}=1$. A gauge transformation (\ref{gauge31}) acts on $\lambda$ as follows:
\begin{equation}
\lambda\rightarrow e^{i(\zeta+\xi-\tau_2\zeta-\tau_1\xi)}\lambda.\,
\end{equation}
hence we can set $\lambda=1$. Therefore, the functions $\varphi_i$ and $\psi_i$ satisfy the system
(\ref{dirac31}). This proves the proposition.

We note that the $\mathbb{R}^{2,1}$-case can be obtained as a reduction by setting $\varphi_2=\bar{\varphi}_1$,
$\psi_2=\bar{\psi}_1$.

\subsection{The $\mathbb{R}^{2,2}$-case}

\begin{proposition}

Suppose that the functions $\varphi_i$, $\psi_i$, $i=1,2$ satisfy the following discrete Dirac
equation:
\begin{equation}
\begin{array}{c}
\tau_2\varphi_1=\alpha\varphi_1+\beta\psi_1,\quad \tau_2\varphi_2=\delta\varphi_2+\gamma\psi_2\\
\tau_1\psi_1=\gamma\varphi_1+\delta\psi_1,\quad \tau_1\psi_2=\beta\varphi_2+\alpha\psi_2,\\
\bar{\alpha}=\alpha,\quad\bar{\beta}=\beta,\quad\bar{\gamma}=\gamma,\quad\bar{\delta}=\delta,\quad
\alpha\delta-\beta\gamma=1.\end{array}
\label{dirac22}
\end{equation}
Then the formulas
\begin{equation}
X_1=\frac{1}{2}\displaystyle\sum\left[(\varphi_1\varphi_2+\bar{\varphi}_1\bar{\varphi}_2)\Delta_1+
(\psi_1\psi_2+\bar{\psi}_1\bar{\psi}_2)\Delta_2\right],
\label{sum22a}
\end{equation}
\begin{equation}
X_2=\frac{i}{2}\displaystyle\sum\left[(\varphi_1\varphi_2-\bar{\varphi}_1\bar{\varphi}_2)\Delta_1+
(\psi_1\psi_2-\bar{\psi}_1\bar{\psi}_2)\Delta_2\right],
\label{sum22b}
\end{equation}
\begin{equation}
X_3=\frac{1}{2}\displaystyle\sum\left[(\varphi_1\bar{\varphi}_2+\bar{\varphi}_1\varphi_2)\Delta_1+
(\psi_1\bar{\psi}_2+\bar{\psi}_1\psi_2)\Delta_2\right],
\label{sum22c}
\end{equation}
\begin{equation}
X_4=\frac{i}{2}\displaystyle\sum\left[(\varphi_1\bar{\varphi}_2-\bar{\varphi}_1\varphi_2)\Delta_1+
(\psi_1\bar{\psi}_2-\bar{\psi}_1\psi_2)\Delta_2\right],
\label{sum22d}
\end{equation}
define an isotropic discrete surface $\vec{X}:\mathbb{Z}^2\rightarrow\mathbb{R}^{2,2}$. Conversely, if
$\vec{X}:\mathbb{Z}^2\rightarrow\mathbb{R}^{2,2}$ is a non-degenerate isotropic discrete surface,
then there exist functions $\varphi_i$, $\psi_i$, $i=1,2$ satisfying equation (\ref{dirac22}) such that equations (\ref{sum22a})-(\ref{sum22d}) hold.
\end{proposition}

\noindent {\bf Proof.} Given functions $\varphi_i$, $\psi_i$, $i=1,2$ satisfying (\ref{dirac22}), a direct calculation
shows that the edges given by equations (\ref{sum22a})-(\ref{sum22d}) are isotropic (\ref{isotropic}) and satisfy
the consistency condition (\ref{consistency}), hence define an isotropic discrete surface in $\mathbb{R}^{2,2}$.

Conversely, suppose that $\vec{X}:\mathbb{Z}^2\rightarrow\mathbb{R}^{2,2}$ is a non-degenerate isotropic discrete surface.
The edges of the lattice satisfy the equations
\begin{equation}
F_1^2+F_2^2=F_3^2+F_4^2,\quad G_1^2+G_2^2=G_3^2+G_4^2,
\end{equation}
therefore, there exist functions $\varphi_i$ and $\psi_i$, where $i=1,2$, such that the edges
are given by the formulas (\ref{sum22a})-(\ref{sum22d}). These functions are defined up to the following local gauge equivalence:
\begin{equation}
\varphi_1\rightarrow \mu\varphi_1,\quad\varphi_2\rightarrow \mu^{-1}\varphi_2,\quad
\psi_1\rightarrow \nu\psi_1,\quad\psi_2\rightarrow \nu^{-1}\psi_2,
\label{gauge22}
\end{equation}
where $\mu$ and $\nu$ are real-valued functions.
The consistency condition (\ref{consistency}) implies that these functions satisfy the following equations
\begin{equation}
(\tau_2\varphi_1)(\tau_2\varphi_2)-\varphi_1\varphi_2=(\tau_1\psi_1)(\tau_1\psi_2)-\psi_1\psi_2,
\label{consistency22a}
\end{equation}
\begin{equation}
(\tau_2\varphi_1)(\tau_2\bar{\varphi}_2)-\varphi_1\bar{\varphi}_2=(\tau_1\psi_1)(\tau_1\bar{\psi}_2)-\psi_1\bar{\psi}_2,
\label{consistency22b}
\end{equation}
and the non-degeneracy condition implies that
\begin{equation}
\varphi_1\bar{\psi}_1-\bar{\varphi}_1\psi_1\neq 0,\quad\varphi_2\bar{\psi}_2-\bar{\varphi}_2\psi_2\neq 0.
\end{equation}

The above equations imply that there exist unique real-valued functions $\alpha_i$, $\beta_i$, $\gamma_i$ and $\delta_i$,
where $i=1,2$, such that the following system of equations is satisfied:
\begin{equation}
\tau_2\varphi_i=\alpha_i\varphi_i+\beta_i\psi_i,\quad t_1\psi_i=\gamma_i\varphi_i+\delta_i\psi_i,\quad i=1,2.
\end{equation}
We can explicitly solve these to obtain
\begin{equation}
\alpha_i=\frac{\bar{\psi}_i(\tau_2\varphi_i)-\psi_i(\tau_2\bar{\varphi}_i)}{\varphi_i\bar{\psi}_i-\bar{\varphi}_i\psi_i},\quad
\beta_i=\frac{\varphi_i(\tau_2\bar{\varphi}_i)-\bar{\varphi}_i(\tau_2\varphi_i)}{\varphi_i\bar{\psi}_i-\bar{\varphi}_i\psi_i},
\end{equation}
\begin{equation}
\gamma_i=\frac{\bar{\psi}_i(\tau_1\psi_i)-\psi_i(\tau_1\bar{\psi}_i)}{\varphi_i\bar{\psi}_i-\bar{\varphi}_i\psi_i},\quad
\delta_i=\frac{\varphi_i(\tau_1\bar{\psi}_i)-\bar{\varphi}_i(\tau_1\psi_i)}{\varphi_i\bar{\psi}_i-\bar{\varphi}_i\psi_i},
\end{equation}
and a direct calculation using (\ref{consistency22a})-(\ref{consistency22b}) shows that
\begin{equation}
\alpha_1\alpha_2-\gamma_1\gamma_2=1,\quad\delta_1\delta_2-\beta_1\beta_2=1,\quad \alpha_1\beta_2-\gamma_1\delta_2=0,
\quad \alpha_2\beta_1-\gamma_2\delta_1=0.
\end{equation}
Solving this system we get that $\alpha_2=\lambda\delta_1$, $\beta_2=\lambda\gamma_1$, $\gamma_2=\lambda\beta_1$
and $\delta_2=\lambda\alpha_1$. A gauge transformation (\ref{gauge22}) acts on $\lambda$ as follows:
\begin{equation}
\lambda\rightarrow (\tau_2\mu)(\tau_1\nu)\mu^{-1}\nu^{-1}\lambda.\,
\end{equation}
hence we can set $\lambda=1$. Therefore, the functions $\varphi_i$ and $\psi_i$ satisfy the system (\ref{dirac22}).
This proves the proposition.

\section{The continuous limit}

\setcounter{equation}{0}

In this section we show that in the continuous limit, the reductions (\ref{dirac21}), (\ref{dirac31}), (\ref{dirac22})
 of the Dirac operator (\ref{gen-discrete-Dirac}) converge to their continuous counterparts.
 
First, consider the operator (\ref{dirac21}). Let $h$ denote the size of the mesh, so that
\begin{equation}
\varphi(x,y+h)=\alpha(x,y)\varphi(x,y)+\beta(x,y)\psi(x,y),\quad
\psi(x+h,y)=\beta(x,y)\varphi(x,y)+\alpha(x,y)\psi(x,y),
\end{equation}
where $\alpha$ and $\beta$ are real and $\alpha^2-\beta^2=1$.
Setting $\beta=hp$, we get that $\alpha=1+O(h^2)$, and expanding the above equation up to $O(h^2)$ gives us
\begin{equation}
\varphi+h\partial_y\varphi=\varphi+hp\psi+O(h^2),\quad
\psi+h\partial_x\psi=hp\varphi+\psi+O(h^2),
\end{equation}
so in the limit $h\rightarrow 0$ we get equation (\ref{cdirac21}).

Similarly, for the operator (\ref{dirac31}) introducing mesh size $h$ we get
\begin{equation}
\varphi(x,y+h)=\alpha(x,y)\varphi(x,y)+\beta(x,y)\psi(x,y),\quad
\psi(x+h,y)=\bar{\beta}(x,y)\varphi(x,y)+\bar{\alpha}(x,y)\psi(x,y),
\end{equation}
where $|\alpha|^2-|\beta|^2=1$. Again, setting $\beta=hp$ gives is $\alpha=1+O(h^2)$, and 
expanding the above equation up to $O(h^2)$ gives us
\begin{equation}
\varphi+h\partial_y\varphi=\varphi+hp\psi+O(h^2),\quad
\psi+h\partial_x\psi=h\bar{p}\varphi+\psi+O(h^2),
\end{equation}
so in the limit $h\rightarrow 0$ we get equation (\ref{cdirac31}).

For the operator (\ref{dirac22}), we introduce a mesh size $h$ to get
\begin{equation}
\varphi_1(x,y+h)=\alpha(x,y)\varphi_1(x,y)+\beta(x,y)\psi_1(x,y),\quad
\psi_1(x+h,y)=\gamma(x,y)\varphi_1(x,y)+\delta(x,y)\psi_1(x,y),
\end{equation}
\begin{equation}
\varphi_2(x,y+h)=\delta(x,y)\varphi_1(x,y)+\gamma(x,y)\psi_1(x,y),\quad
\psi_2(x+h,y)=\beta(x,y)\varphi_1(x,y)+\alpha(x,y)\psi_1(x,y).
\end{equation}
We now use the remaining gauge symmetry (\ref{gauge22}) to set $\alpha=\delta$.
Therefore, if we have a mesh size of $h$, then setting $\beta=hp$, $\gamma=hq$, we see that $\alpha=\delta=1+O(h^2)$, so
expanding the above equation up to $O(h^2)$ gives us 
\begin{equation}
\varphi_1+h\partial_y\varphi_1=\varphi_1+hp\psi_1+O(h^2),\quad
\psi_1+h\partial_x\psi_1=hq\varphi_1+\psi_1+O(h^2),
\end{equation}
\begin{equation}
\varphi_1+h\partial_y\varphi_1=\varphi_1+hq\psi_1+O(h^2),\quad
\psi_1+h\partial_x\psi_1=hp\varphi_1+\psi_1+O(h^2),
\end{equation}
so in the limit $h\rightarrow 0$ we get equation (\ref{cdirac22}).

\section{Acknowledgments}

This work was performed at the Technische Universit\"at in Berlin in June 2009 with financial support from the Polyhedral Surfaces research group. The author would like to thank I.~M.~Krichever for suggesting the problem,  and B.~G.~Konopelchenko and A.~I.~Bobenko for useful discussions.

\end{document}